\documentclass[a4paper,12pt]{article}
\usepackage{amsfonts,amsmath,amsthm}

\def\dR{\mathbb R}
\def\bE{\mathbf E}
\def\bP{\mathbf P}
\def\de{{\rm d}}

\newtheorem{theorem}{Theorem}[section]

\newtheorem{proposition}[theorem]{Proposition}

\begin{document}               
\title{Statistical analysis of the inhomogeneous telegrapher's process}
\author{Stefano Maria Iacus\\Department of Economics, University of 
Milan\\ Via Conservatorio 7, I-20123 Milan - Italy\\ email: 
{\em stefano.iacus@unimi.it}}
\date{}
\maketitle

\begin{abstract}
We consider a problem of estimation for the 
telegrapher's process on the line, say $X(t)$, driven by a Poisson 
process with non constant rate. It turns out that the 
finite-dimensional law of the 
process $X(t)$ is a solution to the telegraph equation with non 
constant coefficients.  We give the explicit law ($\bP_\theta$)
of the process $X(t)$ for a parametric class of intensity functions 
for the Poisson process. We propose an estimator for the parameter 
$\theta$ of $\bP_\theta$
 and we discuss its properties as a first attempt to 
apply statistics to these models.
\end{abstract}



\section{Introduction}
We consider a class of random processes governed by hyperbolic 
equations. These kind of processes have been proposed in the literature 
to describe motions of particles with finite velocities as opposed to 
diffusion-type models. The first contribution in this area goes back 
to Goldstein (1951). He considered the simplest random 
evolution on the real line where a particle, placed in the origin at 
time $0$,  moves with two finite 
velocities $\pm c$ changing its current velocity according to a 
Poisson process of constant rate $\lambda$. He found that the 
distribution of the position of the particle $x$ at time $t$, is a 
solution to the telegraph equation, namely
\begin{equation}
u_{tt}(x,t)+2\lambda u_t(x,t) = c^2 u_{xx}(x,t)
\tag{*}
\label{eq:gold}
\end{equation}
where $u_{tt}(x,t) = \frac{\partial^2}{\partial t^2} u(x,t)$ and so 
forth.
This model has received a great attention in the last decades. Many 
generalizations have been studied and the probabilistic properties of 
these models have been presented. In particular, the finite 
dimensional and the first passage time laws have been presented in 
an explicit form in a series of paper (see Orsingher 1985, 1990 and 1995, 
Foong 1992, Foong and Kanno, 1994).
Different generalizations in the presence of a finite set of 
velocities $c_i$ and Poisson rates $\lambda_i$ (see e.g. Beghin {\it et 
al} 1999)  as well as motions in two or more dimensions 
(Orsingher 1986 and 2000, Orsingher and Kolesnik 1996, Kolesnik and 
Turbin 1991, Iacus 1995) 
have also been considered. In general, for these models it is quite hard to 
find explicitly the distribution laws of interest.

In this paper we present a generalization of the model when the 
underlying Poisson process is not homogeneous that are interesting 
from a statistical point of view. Up to our knowledge, this is the 
first attempt to apply statistics in this field.

The paper is divided into two parts: the first concerns the 
probabilistic analysis of the model and the latter is devoted to 
statistical estimation.  In Section \ref{sec:vel} we introduce the 
model, we study the so called velocity process and we derive the 
following partial 
differential equation analogous to \eqref{eq:gold} :
\begin{equation}
u_{tt}(x,t)+2\lambda(t) u_t(x,t) = c^2 u_{xx}(x,t)
\tag{**}
\label{eq:mia}
\end{equation}
where $\lambda(\cdot)$ is the intensity function of the Poisson process 
(see e.g. Kutoyants 1998).
Then (in \S \ref{sec:model}) we give conditions on $\lambda(\cdot)$ under 
which an explicit solution to  \eqref{eq:mia} can be found. It turns 
out that it is a parametric family of intensity measures
$\lambda(\cdot)=\lambda_\theta(\cdot)$. In the second part of 
the paper (Section \ref{sec:est}) we consider the problem of estimation 
of the parameter $\theta$ for this particular class of solutions and 
we present asymptotically efficient --in the minimax sense-- estimators for $\theta$.

\section{The model and the velocity process}\label{sec:vel}
Consider a particle placed in $x=0$  at time $t=0$. It can move leftward 
or rightward with finite velocity $c$. Changes of direction occur at 
Poisson times. This means that the particle will move at constant 
speed (either $+c$ or $-c$) among two successive Poisson events.
We introduce the velocity process $V(t)$: 
\begin{equation}
V(t)= V(0)(-1)^{N(t)},\quad t>0\,,
\label{eq:vel}
\end{equation}
where $V(0)$ is  a random variable taking values $\pm c$ with probability 
$\frac12$ and independent from the Poisson process
$N(t)$ with intensity measure $\Lambda(\cdot)$ :
$$
\Lambda(t) = \int_0^t \lambda(s) \de s\,,\qquad t>0\,,
$$
where $\lambda(\cdot)$ -- the intensity function -- belongs to $\mathcal C^1(\mathbb R)$. 
The increments $N(t)-N(s)$, $0<s<t$, are then 
distributed according to the Poisson law with parameter
$\Lambda(t)-\Lambda(s)$. 

Our main goal is to determine the law governing the process $X(t)$
that represents the position of the particle at time $t$ : 
\begin{equation}
X(t)=V(0)\int_0^t (-1)^{N(s)}\de s\,.
\label{processo}
\end{equation}
We start with the analysis of the velocity process $V(t)$, $t>0$. At 
first, we note that
the two probabilities
$$
\begin{aligned}
p^{(c)}(t) &= \bP(V(t)=c)    \\
p^{(-c)}(t) &= \bP(V(t)=-c)
\end{aligned}   
$$
are solutions of the following  system of differential equations
\begin{equation}
\left\{\begin{aligned}
p^{(c)}_t(t) &= \Lambda_t(t)(p^{(-c)}(t) - p^{(c)}(t))\\
p^{(-c)}_t(t) &= \Lambda_t(t)(p^{(c)}(t) - p^{(-c)}(t))=-p^{(c)}_t(t)
\end{aligned}   
\right.
\label{sis1}
\end{equation}
This can be proved by Taylor expansion of the two functions $p^{(c)}_t(t)$ and 
$p^{(-c)}_t(t)$.
The following conditional laws: 
$$
\bP(V(t)=c|V(0)=c) \qquad\text{and}\qquad \bP(V(t)=-c|V(0)=c) \,,
$$
characterize the velocity process. Their explicit form are as follows.
\begin{proposition}
\begin{equation}
\bP(V(t)=c|V(0)=c) = 1 - \int_0^t \lambda(s) e^{-2\Lambda(s)} \de s =
\frac{1+e^{-2\Lambda(t)}}{2}\,,
\label{eq:vel1}
\end{equation}
\begin{equation}
\bP(V(t)=-c|V(0)=c) =  \int_0^t \lambda(s) e^{-2\Lambda(s)} \de s =
\frac{1-e^{-2\Lambda(t)}}{2}
\label{eq:vel2}
\end{equation}
\end{proposition}
\begin{proof}
We give only the derivation of \eqref{eq:vel1} because the other 
follows by symmetry. Conditioning on $V(0)=c$ implies that:
 $p^{(c)}(0) = 1$, $p^{(c)}_t(0)=-\Lambda_t(0)=-\lambda(0)$, 
 $p^{(-c)}(0)=0$ and  $p^{(-c)}_t(0)=\lambda(0)$.
From \eqref{sis1} we obtain that
$$
\frac{p_{tt}^{(c)}(t)}{p_t^{(c)}(t)}=\frac{\de }{\de t} \log p^{(c)}_t(t) =-\frac{2 
\Lambda_t(t)^2-\Lambda_{tt}(t)}{\Lambda(t)}\,,
$$
and by simple integration by parts it emerges that:
$$
p^{(c)}_t(t) = \exp\left\{
- 2 \Lambda(t) + \log \left(\frac{\lambda(t)}{\lambda(0)}\right) + 
\log (-\lambda(0))
\right\} \,,
$$
from which follows \eqref{eq:vel1}.
\end{proof}
Remark that it is also possible to derive  \eqref{eq:vel1} by simply 
noting that
$$
\begin{aligned}
\bP(V(t)=c|V(0)=c) &= \bP\left(
\bigcup_{k=0}^\infty N(t) =2k
\right)\\
&=\sum_{k=0}^\infty  \frac{\Lambda(t)^{2k}}{(2k)!} e^{-\Lambda(t)}\\
&=\frac{1+e^{-2\Lambda(t)}}{2}
\end{aligned}
$$
and then $\bP(V(t)=-c|V(0)=c) = 1 -\bP(V(t)=c|V(0)=c)$.

We conclude the analysis of the velocity process by giving also the 
covariance function of the couple $(V(t),V(s))$, $s,t>0$. In fact, it 
can be easily proven that the characteristic function of 
the couple $(V(t),V(s))$ is, for all  $(\alpha,\beta)\in\dR^2$,
$$
\bE\left(e^{i\alpha V(s)+i\beta V(t)}\right)=\cos(\alpha 
c)\cos(\beta c)-e^{-2(\Lambda(t)-\Lambda(s))}
\sin(\alpha c)\sin(\beta c),
$$
and the covariance function is then
$$
\bE\left(V(s) V(t)\right) = \left.
-\frac{\partial^2}{\partial\alpha\partial\beta}
\bE\left(e^{i\alpha V(s)+i\alpha 
V(t)}\right)\right|_{\alpha=\beta=0}=c^2 
e^{-2|\Lambda(t)-\Lambda(s)|}\,.
$$

\subsection{Derivation of the telegraph equation}
In order to analyze the distribution of the position of the 
particle 
\begin{equation}
P(x,t) = \bP(X(t)<x)   
\label{eq:P}
\end{equation}
we introduce the two distribution functions
$
F(x,t) = \bP\left(X(t) < x , V(t) = c\right)$ and
$B(x,t) = \bP\left(X(t) < x , V(t) = -c\right)$, so that
$
P(x,t) = F(x,t) + B(x,t)$ and $W(x,t) = F(x,t) - B(x,t)$. 
The function $W(\cdot,\cdot)$ is usually called the ``flow function".
Next result gives the analogous to \eqref{eq:gold} in the case of 
nonhomogeneous Poisson process.
\begin{proposition} Suppose that $F(\cdot,\cdot)$ and $B(\cdot,\cdot)$ are two times 
differentiable in $x$ and $t$, then
\begin{equation}
\left\{
\begin{aligned}
F_t(x,t) &= -c F_x(x,t) -\lambda(t)(F(x,t)-B(x,t))\\
B_t(x,t) &= c B_x(x,t) + \lambda(t)(F(x,t)-B(x,t))
\end{aligned}   
\right.
\label{sis:fb}
\end{equation}
moreover $P(x,t)$ is a solution to the following telegraph equation with non constant coefficients
\begin{equation}
\frac{\partial^2}{\partial t^2} u(x,t) + 2 \lambda(t) 
\frac{\partial}{\partial t} u(x,t) = c^2 \frac{\partial^2}{\partial x^2} 
u(x,t)\,.
\label{eq:tel}
\end{equation}
\end{proposition}
\begin{proof}
By Taylor expansion, one gets that $F(\cdot,\cdot)$ and $B(\cdot,\cdot)$ are 
solutions to \eqref{sis:fb} and
rewriting system \eqref{sis:fb} in terms of the functions 
$W(\cdot,\cdot)$ 
and $P(\cdot,\cdot)$ it emerges that
\begin{equation}
\left\{\begin{aligned}
P_t(x,t) &= -c W_x(x,t) \\
W_t(x,t) &= c P_x(x,t) - 2 \lambda(t) W(x,t)
\end{aligned}   
\right.
\label{sis:pw}
\end{equation}
The conclusion arises  by direct substitutions.
In fact, by the first of system \eqref{sis:pw} we have
$$
P_{tt}(x,t) = \frac{\partial}{\partial t} P_t(x,t)
 =  -c \frac{\partial}{\partial t}  W_x(x,t)\,.
$$
Furthermore
$$
\begin{aligned}
W_{tx}(x,t) &= \frac{\partial}{\partial x} (c P_x(x,t) - 2 \lambda(t) 
W(x,t)) \\
&= c P_{xx}(x,t) - 2\lambda(t) W_x(x,t)\\
&= c P_{xx}(x,t) +\frac{2}{c} \lambda(t) P_t(x,t)
\end{aligned}
$$
by using respectively the second and the first equation of system \eqref{sis:pw}.
\end{proof}

\subsection{The explicit law of the telegraph process}\label{sec:model}
In this section we give the explicit form of distribution function 
\eqref{eq:P} for a particular class of intensity functions. The idea is to reduce 
equation \eqref{eq:tel} to  a partial differential equation for which the 
general solution is available. This is done by imposing conditions on the family of intensity functions 
$\lambda(\cdot)$. Here we give only one type of solution that is interesting 
from the statistical point of view. The result presented in the next theorem is interesting in itself.
\begin{theorem}
Suppose that the intensity function of the Poisson process $N(t)$ in 
\eqref{processo} is
\begin{equation}
\lambda(t) = \lambda_\theta(t) = \theta \tanh(\theta t)
,\quad \theta\in\mathbb R\,.
\label{eq:lambda}
\end{equation}
Then, the absolutely continuous component $p_\theta(\cdot)$  of distribution \eqref{eq:P}, 
 is given by
\begin{equation}
p_\theta(x,t,c) =
\begin{cases}
\frac{\theta \, t}{\cosh(\theta \, t)}    
\frac{I_1\left(\frac{\theta}{c} \sqrt{{c^2}\, 
{t^2}-{x^2}}\, \right)}{2{\sqrt{{c^2}\, {t^2}-{x^2}}}},&|x|<ct\\
0,& \text{otherwise}
\end{cases}
\end{equation}
\end{theorem}    
\begin{proof}
We start without assuming \eqref{eq:lambda} and by noting that  $p(x,t)$ is a solution 
to \eqref{eq:tel}. Moreover, 
the process $X(t)\in(-ct,ct)$ only if there occur at least one 
Poisson event up to time $t$, thus
\begin{equation}
\int_{-ct}^{ct} p(x,t) \de x = 1 - \bP(N(t)=0) = 1 - e^{-\Lambda(t)}
\label{eq:bound}
\end{equation}
that gives one of the conditions to solve \eqref{eq:tel}.
We now search for solutions of the following type:
$$
v(x,t) = e^{\Lambda(t)} u(x,t)\,.
$$
Thus the function $v(\cdot,\cdot)$ satisfies the following partial differential equation
\begin{equation}
v_{tt}(x,t) - v(x,t)\left (\lambda'(t) + \lambda^2(t)\right) = c^2 
v_{xx}(x,t)\,.
\label{eq:para}
\end{equation}
For a generic function $\lambda(\cdot)$ a solution to 
\eqref{eq:para} 
not available. 
If the intensity function $\lambda(\cdot)$ 
satisfies the following ordinary differential equation
\begin{equation}
\lambda'(t) + \lambda^2(t) = \theta^2,\quad 
t>0,\quad\theta\in\mathbb R\,,
\label{eq:cool}
\end{equation}
then equation \eqref{eq:para} becomes
\begin{equation}
v_{tt}(x,t) - v(x,t) \theta^2 = c^2 v_{xx}(x,t)
\label{eq:ors}
\end{equation}
and the solution to it can be written 
explicitly\footnote{Remark that, if one solves the 
problem $\lambda'(t) + \lambda^2(t) = 0$
equation \eqref{eq:para} reduces to the standard one-dimensional wave 
equation $v_{xx}=c^2 v_{xx}$ and its solution can also be determined 
but it is of no statistical interest.}. 
With the initial condition $\lambda(0) = 0$ the solution to equation  \eqref{eq:cool} is condition
\eqref{eq:lambda} and then $\Lambda_\theta(t) = \ln(\cosh(\theta t))$.
Following Orsingher (1985), by the change of variable $s=\sqrt{c^2 t^2 - 
x^2}$ we transform equation \eqref{eq:ors} into the following standard Bessel's 
equation
\begin{equation}
v_{ss}+\frac{1}{s}v_s - \left(\frac{\theta}{c}\right)^2 v = 0
\label{eq:bsl}
\end{equation}
whose general integral is
$$
v(x,t)=A I_0\left(
\frac{\theta}{c}\sqrt{c^2 t^2 - x ^2}
\right)+ B K_0\left(
\frac{\theta}{c}\sqrt{c^2 t^2 - x ^2}
\right)
,\quad |x|<ct\,.
$$
The function $I_k(x)$ is the modified Bessel function of first kind and order k
while $K_0(x)$ is the second type Bessel function of order 0 with the 
unpleasant property that $\lim_{x\to 0^+} K_0(x)=\infty$, so we put 
$B=0$. In terms of $p_\theta(\cdot,\cdot)$ 
the solution is of the following form
\begin{equation}
p_\theta(x,t) = K e^{-\Lambda_\theta(t)} I_0\left(
\frac{\theta}{c}\sqrt{c^2 t^2 - x ^2}
\right),\quad |x|<ct\,. 
\label{eq:boh}
\end{equation}
From \eqref{eq:bound} it follows that
$$
\int\limits_{-ct}^{ct} p_\theta(x,t) \de x = 1 - e^{-\Lambda_\theta(t)}  
= 1 - \frac{1}{\cosh(\theta t)} = 1 - \frac{2}{e^{^-\theta t}+e^{\theta t}}
$$
and so, there is no constant value $K$  that satisfies \eqref{eq:boh}.
To turn around this problem, we note that if $I_0\left(
\frac{\theta}{c}\sqrt{c^2 t^2 - x ^2}
\right)$ is a solution to \eqref{eq:bsl} so is its partial derivative 
with respect to $t$. Hence we search for solutions of the form
$$
p_\theta(x,t) = e^{-\Lambda_\theta(t)} \left(
A I_0\left(
\frac{\theta}{c}\sqrt{c^2 t^2 - x ^2}
\right) + B \frac{\partial}{\partial t}I_0\left(
\frac{\theta}{c}\sqrt{c^2 t^2 - x ^2}
\right)
\right)\,.
$$
Form the following two equalities (see Orsingher, 1995) 
$$
\int_{-ct}^{ct} I_0\left(
\frac{\theta}{c}\sqrt{c^2 t^2 - x ^2}
\right)\de x=\frac{\theta}{c}\left(e^{\theta t}-e^{-\theta t}\right)
$$
and
$$
\int_{-ct}^{ct} 
\frac{\partial}{\partial t} I_0\left(
\frac{\theta}{c}\sqrt{c^2 t^2 - x ^2}
\right)\de x=c\left(e^{\theta t}+e^{-\theta t}\right)-2c
$$
condition \eqref{eq:bound} implies that
$$
1 - e^{-\Lambda_\theta(t)} = 
e^{-\Lambda_\theta(t)} \left(
A (1-\lambda_\theta(t)) \frac{c}{\theta} 
\left(e^{\theta t}-e^{-\theta t}\right) + B c 
\left(e^{\theta t}+e^{-\theta t}-2\right)
\right)\,.
$$
Observe that for $\lambda_\theta(t)$ as in \eqref{eq:lambda} we have
$
1-\lambda_\theta(t) = 1  - \theta\,\frac{e^{\theta t}-e^{-\theta t}}{e^{-\theta 
t}+e^{\theta t}}
$
thus the equality in the above equation is attained by taking $A=0$ and $B=\frac{1}{2c}$.
The solution is finally
$$
p_\theta(x,t) =
\frac{\theta \, t \,I_1\left(\frac{\theta}{c} \sqrt{{c^2}\, 
{t^2}-{x^2}}\, \right)}{({e^{-\theta\,t }}+{e^{\theta\,t 
}}){\sqrt{{c^2}\, {t^2}-{x^2}}}},\quad |x|<ct\\
$$
that is non negative and satisfies the conditions required. It also 
verifies equation \eqref{eq:tel} and this is a boring calculus' exercise.
\end{proof}
Note that the non-null component of \eqref{eq:P} can be written as
$$\frac{1}{2 \,\cosh(\theta \, t)}
\left( \left(\delta(x-ct)+\delta(x+ct)\right)
+ \frac{\theta \, t I_1\left(\frac{\theta}{c} \sqrt{{c^2}\, 
{t^2}-{x^2}}\, \right)}{{\sqrt{{c^2}\, {t^2}-{x^2}}}}\right),\,\, 
x\in[-ct,ct]
$$
where $\delta(x)$ is the Dirac delta.
By the result of last theorem we can now state for the seek of 
completeness the following result.
\begin{proposition}
Suppose that the intensity function of the Poisson process $N(t)$ in 
\eqref{processo} is
$$\lambda(t) = \lambda_\theta(t) = \theta \tanh(\theta t)
,\quad \theta\in\mathbb R\,.
$$
Then, the distribution function \eqref{eq:P}  is 
$$
P_\theta(x,t)=
\begin{cases}
0,& x<-ct\\
\frac{1}{2\,\cosh(\theta \, t)} + \int\limits_{-ct}^x p_\theta(u,t)\de 
u,& -ct\leq x<ct\\
1,&x\geq ct
\end{cases}    
$$
\end{proposition}    
\section{Parameter estimation}\label{sec:est}
We consider the problem of estimation of the parameter $\theta$ 
for the model discussed in \S\ref{sec:model}.  We suppose that the velocity $c$ is known, 
if it is not, 
one can easily determine without error its value, by observing two 
successive switches of velocity (or just the first switch) and the time 
between the two occurrences.
Two possibile schemes of observations can be considered:  {\it i}) one trajectory is observed up to some time $T$, letting 
$T\to\infty$, and {\it ii}) $n$ 
independent and identically distributed observation $X_i(t)$ of the 
trajectories are observed on a fixed time interval $[0,T]$, letting $n\to\infty$.
Instead of working on the law of $X(t)$ we do inference on $\theta$ 
via the Poisson process.

We use the  method of moments to 
estimate $\theta$.
For scheme {\it i}), recall that $\bE N(t)=\Lambda_\theta(t)$
thus, if $\pi_T$ is the number of observed switches up to time $T$, it 
suffices to find the solution to
$$
\Lambda_\theta(T)  = \ln({\rm cosh}(\theta T)) = \pi_T
$$
that is
$$
\tilde\theta_T = \theta(\pi_T) = \frac1{T}\,{\rm 
arcosh}\left(e^{\pi_T}\right)\,.
$$
Good asymptotic properties of this estimator cannot be established without 
assuming more conditions on the Poisson process. We give now the 
optimal solution in terms of the second scheme.

Suppose now to be able to observe $n$ independent copies $X_i(t)$ of the process 
$X(t)$ up to a fixed horizon $T$, as to say, we observe the trajectories of $n$ 
particles that do not interact. Denote by $\pi_i$ the number of switches in 
each replication. Then 
$$
\hat{\pi}_n = \frac1{n} \sum_{i=1}^n \pi_i
$$
is a consistent estimator of $\pi = \Lambda_\theta(T)$ as $n\to\infty$, 
moreover $\sqrt{n}(\hat{\pi}_n - \pi)$ is asymptotically 
Gaussian $\mathcal N(0,\pi)$. The estimator 
$\hat{\pi}_n$ -- called empirical measure -- is also 
asymptotically efficient in the minimax sense (for these and many other 
results see Kutoyants, 1998). In this particular case, an asymptotically 
efficient estimator for $\theta$ exists. In 
fact,
$$
\pi = \ln({\rm cosh}(\theta T))
\quad\text{and}\quad 
\theta = \theta(\pi)= \frac1{T}\,{\rm arcosh}\left(e^\pi \right)
$$
so the estimator of moments is $\hat\theta_n = \theta(\hat\pi_n)$. By 
the so-called $\delta$-method we desume the asymptotic 
behavior of $\theta(\hat\pi_n)$. In fact
$$
\sqrt{n}(\theta(\hat\pi_n)-\theta(\pi)) = \sqrt{n}(\hat{\pi}_n - \pi) 
\theta'(\pi)+o(|\hat{\pi}_n - \pi|)
$$
where the last term converges to 0 in probability by consistency.
Thus, $\sqrt{n}(\theta(\hat\pi_n)-\theta(\pi))$ has, asymptotically, a  
centered Gaussian law with variance given by 
$$\pi/\theta'(\pi)^2
= \frac{2\, T^2 \pi}{1 + \coth(\pi)} 
$$
or better
$$
\sqrt{n}(\hat\theta_n-\theta) \Longrightarrow \mathcal N\left(0,T^2 
\frac{\ln(\cosh(\theta T))}{\coth(\theta T)^2} \right)
$$
This estimator is also the maximum likelihood estimator for $\theta$ 
that is asymptotically efficient in this regular problem.
To conclude, the estimator $\hat\theta_n$ is consistent for $\theta$ 
and asymptotically Gaussian, moreover it is asymptotically efficient 
for $\theta$ with asymptotic variance given by the inverse of 
the Fisher information in this problem.

\section*{References}
\begin{itemize}
\item Beghin, L., Nieddu, L., Orsingher, E. (1999), Probabilistic 
analysis of the telegrapher's process with drift by means of 
relativistic transformations, to appear in {\em Jour. Appl. Math. and Stoch. 
Anal.}
\item Foong, S.K. (1992), First passage time, maximum displacement and 
Kac's solution of the telegraphers equation, {\em Physical Review A}, 
{\bf 46}, 707-710.
\item Foong, S.K., Kanno S. (1994), Properties of the telegrapher's random 
process with or without a trap, {\em Stoch. Proc. and Their Appl.}, 147-173.
\item Goldstein, S. (1951), On diffusion by discontinuous movements 
and the telegraph equation, {\em Quart. J. Mech. Appl. Math.}, {\bf 4}, 129-156.
\item Iacus, S.M. (1995), Random motions governed by hyperbolic 
equations, in  {\em Acts of European Young Statisticians Meeting EYSM95}, 
(Eramus University, Rotterdam),  44-49.
\item Kolesnik, A.D., Turbin, A.F. (1991), Infinitesimal hyperbolic 
operator of the Markovian random evolutions in $R^n$, {\em Doklad. Akadem. 
Nauk. Ukrain}, 11-14 (in russian).
\item Kutoyants, Yu.A. (1998), {\em Statistical inference for spatial 
Poisson pro\-ces\-ses}, Lecture Notes in Statistics (Springer-Verlag).
\item Orsingher, E. (1985), Hyperbolic equations arising in random 
models, {\em Stoch. Proc. and Their Appl.}, 49-66.
\item Orsingher, E. (1986), A planar random motion governed by the 
two-dimensional telegraph equation, {\em Jour. Appl. Prob.}, {\bf 22}, 385-397.
\item Orsingher, E. (1990), Probability law, flow function, maximum 
distribution of wave-governed random motions and their connections 
with Kirchoff's laws, {\em Stoch. Proc. and Their Appl.}, 49-66.
\item Orsingher, E. (1995), Motions with reflecting and absorbing 
barriers driven by the telegraph equation, {\em Random Operators and 
Stochastic Equations}, {\bf 3:1}, 9-21.
\item Orsingher, E. (2000), Exact joint distribution in a model of 
planar random motion,  {\em Stoch. and Stoch. Reports}, {\bf 1-2}, 1-10.
\item Orsingher, E., Kolesnik, A. (1996), The exact probability law of a 
planar random motion governed by a fourth-order hyperbolic equation, 
{\em Theory of Probability and its Appl.}, {\bf 41}, 379-386.
\end{itemize}
\end{document}